# SCOR Paper

**Oskar Laverny**[1,2*]
**Alessandro Ferriero**[2]
**Ecaterina Nisipasu**[2]

# Parametric divisibility of stochastic losses

[1]Insitut Camille Jordan, UMR 5208,
Université Claude Bernard, Lyon 1, France
[2]SCOR SE, France
*Corresponding author: oskar.laverny@gmail.com

SCOR
The Art & Science of Risk

# Abstract

A probability distribution is n-divisible if its $n^{th}$ convolution root exists. While modeling the dependence structure between several (re)insurance losses by an additive risk factor model, the infinite divisibility, that is the n-divisibility for all $n \in \mathbb{N}$, is a very desirable property. Moreover, the capacity to compute the distribution of a piece (i.e., a convolution root) is also desirable. Unfortunately, if many useful distributions are infinitely divisible, computing the distributions of their pieces is usually a challenging task that requires heavy numerical computations. However, in a few selected cases, particularly the Gamma case, the extraction of the distribution of the pieces can be performed fully parametrically, that is with negligible numerical cost and zero error. We show how this neat property of Gamma distributions can be leveraged to approximate the pieces of other distributions, and we provide several illustrations of the resulting algorithms.

# Contents

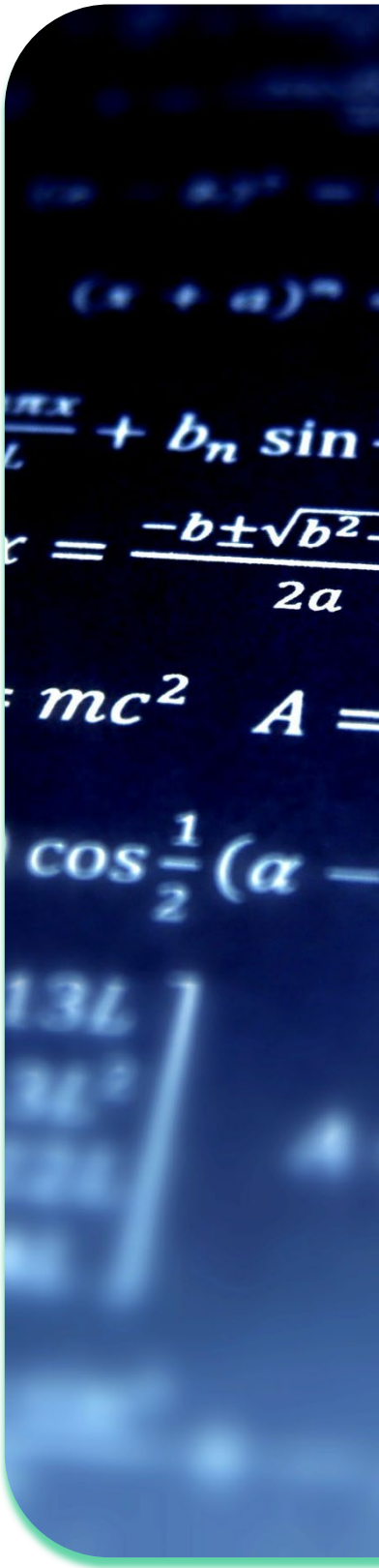

# 1. Introduction

Insurance and reinsurance risks are classically modeled via positive random variables representing loss amounts, whose distributions are estimated from empirical data and/or specific information about the underlying process producing these losses (such as contract details or physical phenomena). Examples of the wide range of literature on modeling methods and associated tools can be found in the Reference section [1–5] at the end of this document. We consider here the internal modeling point of view, where the distributions of several losses $X_1, \ldots, X_d$ are supposed to be known, but the dependence structure between them must still be evaluated and taken into account to assess the variability and the extremal behavior of the total loss $\sum_{i=1}^{d} X_i$.

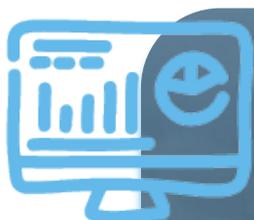

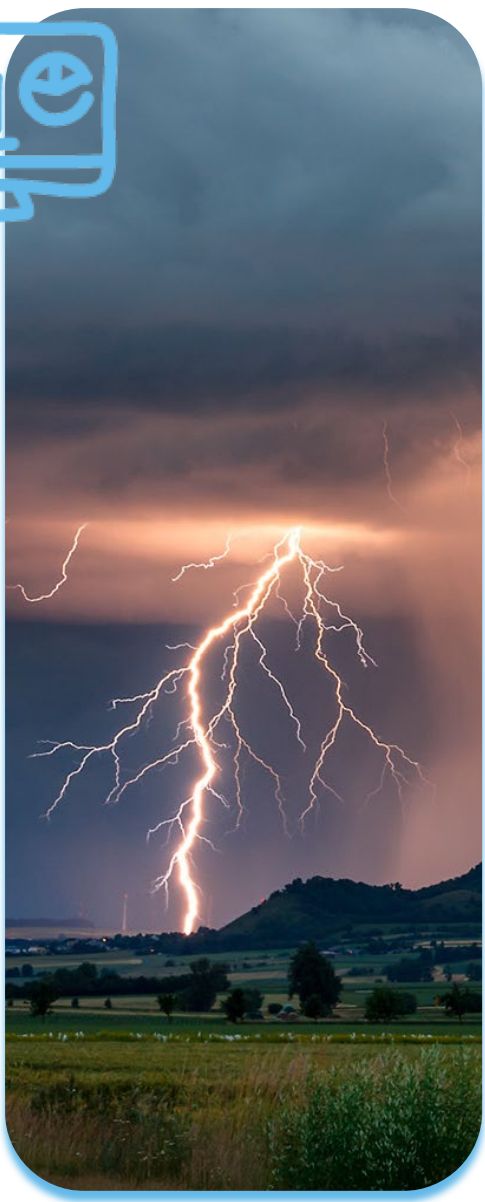

Setting such a dependency is an important but complex matter. Important, since the diversification effect between the $X_i$'s and their potential tail dependencies might induce drastically different behaviors for the total risk, especially in the extremes. Complex, since the quality and quantity of available information is usually not very good.

When there is relevant data about the joint behavior of the marginals, we rely on parametric estimations of copulas. However, when the quantity and/or quality of available data is not sufficient to produce relevant estimations of copulas, we must take the viewpoints of experts into account. Experts' viewpoints might not be in a format adapted to deducing the parametrization of a given copula model, since the interpretability of a Clayton's $\theta$, for example, is hard to grasp. Let alone choosing between several parametric families.

To handle and aggregate several of these experts' viewpoints, we consider the use of an additive risk factor model. These models manage the dependencies by using latent risk factors, which are allowed to produce losses in each of the marginals. More formally, consider that there exist uniform random variables $U_1, \ldots, U_n$, all independent of each other, such that the random variables $X_1, \ldots X_d$ are each written as:

$$X_i = \sum_{j=1}^{n} Q_{i,j}(U_j). \tag{1}$$

Of course, the parameter $n$ is very important to the model: the cases $n \ll d$ or $n \gg d$ are evidently different regimes (the first one inducing the irregularity of the produced random vector **X**, the second one allowing, for example, full independence).

Our notation suggests the fact that $Q_{i,j}$'s are quantile functions, although this is not even necessary for the model to make sense. Some $Q_{i,j}$ might be identically zero, and we tend to hope that many of them will be for sparsity reasons. Note that if all $Q_{i,j}$ are non-decreasing (or if they are all non-increasing), the $j^{th}$ risk factors produce losses $\{Q_{1,j}(U_j), \ldots, Q_{d,j}(U_j)\}$ that are comonotone. For all these reasons, and for the interpretability that arises from this comonotony, we make the assumption that all $Q_{i,j}$ are quantile functions (and thus non-decreasing).

This model, although appealing, is clearly intractable when nothing more is supposed about the set $(Q_{i,j})_{i,j}$. The first simplification, discussed in Section 2, is the restriction to convolutional roots or *pieces*, leveraging the infinite divisibility of the marginals. Another useful restriction is on the class of functions allowed for $Q_{i,j}$'s: an important case is the Gamma restriction, which we detail and extend in Section 3. Section 4 discusses a few examples of applications, and Section 5 concludes.

## 2. Parametric divisibility

Consider again Equation (1), and assume from now on that all $Q_{i,j}$ are quantile functions. Denote the resulting random variables by:

$$Y_{i,j} = Q_{i,j}(U_j). \qquad (2)$$

For any given index $i$, the set of random variables $Y_{i,1}, \ldots, Y_{i,n}$ is clearly a set of independent random variables, since $U_1, \ldots, U_n$ are independent. The first assumption we will make to simplify the model is that, for all $i$, $Y_{i,1}, \ldots, Y_{i,n}$ are all convolutional roots (also called pieces) of $X_i$, according to Definition 1.

**Definition 1** ($n^{th}$-root and $\beta$-piece). *Let $X$ be a positive, continuous, random variable. The $n^{th}$ convolutional root of $X$ (also called the $\frac{1}{n}$-piece of $X$), when it exists, is the common distributions of independent random variables $Y_1, \ldots Y_n$ such that $X = \sum_{i=1}^{n} Y_i$. More generally, if $M$ denotes the moment generating function of $X$, that is $M(t) = \mathbb{E}(e^{tX})$, we define (when it exists) the $\beta$-piece of a random variable $X$ as the distribution having moment generating function $t \mapsto M(t)^{\beta}$.*

With a slight abuse of terminology, we make no distinction between the distribution of the random variable and the random variable itself when this is clear from the context.

When the $n^{th}$-convolutional root exists, we say that $X$ is $n$-divisible. When any $\beta$-piece exists for all $\beta \in ]0,1]$, we say that $X$ is infinitely divisible. By convention, the 0-piece is the constant random variable with value 0. Again, we might speak about the divisibility of a random variable or of its distribution.

The moment generating function of a positive and infinite divisible random variable $X$ can be represented canonically by

$$M(t) = \mathbb{E}(e^{tX}) = \exp\left\{at + \int_{\mathbb{R}_+} (e^{ty} - 1)\, L(dy)\right\},$$

where, according to [6], $a \geq 0$ and $L$, the *Lévy measure*, is non-negative and satisfies $\int_{\mathbb{R}_+} \min(1, y) L(dy) \leq \infty$. The parameter $a$ is the left extremity of the support of the distribution, and the Lévy measure $L$ is uniquely determined. This is essentially a special case of the Lévy-Kitchnine representation that holds for any real infinite divisible distributions, already discussed in [7,8].

The divisibility of a distribution is sometimes a desirable property, especially in (re)insurance risk modeling, when we consider the distribution as a candidate model for aggregated losses. Indeed, standard collective risk models in insurance consider that, when we have $n$ insurers generating independent and identically distributed losses $Y_1, \ldots, Y_n$, the global loss is given as the convolution of these losses,

$$X = \sum_{i=1}^{n} Y_i.$$

Moreover, each insurer may have several losses and therefore the $Y_i$'s could also be expressed as sums of (more or less) independent random variables.

Therefore, due to the very nature of the phenomena we model, insurance losses, the divisibility of the modeling distribution for $X$ is already in our hypothesis set. For this reason, when the time comes to model the dependence structure between several random variables $X_1, \ldots X_d$, corresponding, for example, to different lines of business, using the divisibility property seems natural and appealing to the practitioner.

Indeed, under Equations (1) and (2), the dependence structure of the random vector $\mathbf{X} = (X_1, \ldots X_d)$ can be fully described by the knowledge of all the $Q_{i,j}$ functions, that is all the distributions of $Y_{i,j}$'s. Typically, Ferriero [7,8] proposes to parametrize the problem by a set of constants $\beta_{i,j}$ such that the following assumption holds for all $i \in 1, \ldots, d$:

$$Y_{i,j} \text{ is the } \beta_{i,j}\text{-piece of } X_i. \tag{3}$$

This directly implies that all $\beta_{i,j}$ are in $[0,1]$ (with $\beta_{i,j} = 0$ inducing $Y_{i,j} \equiv 0$), and that

$$\forall i \in 1, \ldots, d, \quad \sum_{j=1}^{n} \beta_{i,j} = 1.$$

The class of available dependence structures in this model is still wide, since the number $n$ of pieces can always be increased, and pieces of different sizes can be constructed. The model contains many interesting behaviors for the final random vector $\mathbf{X}$: asymmetric dependence structures, tail dependencies, conditional independence, etc. See [7] for a broad review of potential models covered by this setting.

A prerequisite of this dependency construct is that we already have known distributions for the losses. Due to the real nature of objects we are modeling, and due to previous considerations, it makes sense to choose for $X_1, \ldots X_d$ distributions that are infinitely divisible. This is the line of reasoning taken by Thorin in 1977.

Thorin was an actuary who noted that common practice, especially in reinsurance, often modeled losses with heavy tailed distributions such as the log-Normal and the Pareto distributions. But these distributions were not known at the time to be infinitely divisible, or even divisible at all. Proving such a fact was a hard problem, which Thorin tackled in two founding articles [9,10], using a purposely made-up class of distributions, the so-called *generalized Gamma convolutions*, which are (briefly) defined as weak limits of convolutions of Gamma distributions. We will see that this class plays a particular role inside the (wider) class of infinitely divisible positive random variables.

Consider again a positive and infinitely divisible random variable $X$. Recall that the $\beta$-piece of $X$ is characterized as follows:

**Remark 2** (Characterization of the $\beta$-piece). *For $X$ a positive and infinitely divisible random variable with moment generating function $M(t)$, for $\beta \in [0,1]$, the $\beta$-piece of $X$ has moment generating function $M_\beta(t)$ such that*

$$M(t) = M_\beta(t)^\beta. \tag{4}$$

Equation (4) must hold for all $t \in \mathbb{C}$ where $M(t)$ is defined. The construction of a piece is therefore closely related to deconvolutional problems. Since $M$ is essentially a complex function, and since $\beta \in [0,1]$, the first branch of the $\beta$-power of $M(t)$ can be taken.

For our purposes of modeling insurance risks with an additive risk factor model based on infinite divisibility, it would be interesting to identify distributions whose pieces are easily identifiable through a parametric model. Indeed, one of the most important operations we might need to do is sample these pieces, and the cost of this sampling must be controlled as much as possible. The sampling algorithm is described in Algorithm 1.

---

**Algorithm 1:** Sampling from a risk factor model.

**Input:** Specification of a risk factor model for $\boldsymbol{X}$, number of samples $N \in \mathbb{N}$.
**Result:** An $N$-sample $\boldsymbol{x} = (\boldsymbol{x}_1, \ldots \boldsymbol{x}_N) \in \mathbb{R}_+^{N \times d}$
**foreach** $j \in 1, \ldots, n$ **do**
  Sample $u_{j,1}, \ldots, u_{j,N}$ from a random uniform.
  **foreach** $i \in 1, \ldots, d$ **do**
    Compute $Q_{i,j}$, the quantile function of the $\beta_{i,j}$-piece of $X_i$.
    **foreach** $k \in 1, \ldots, N$ **do**
      Set $x_{i,k} \mathrel{+}= Q_{i,j}(u_{i,k})$
    **end**
  **end**
**end**
Return $\boldsymbol{x}$.

---

Let us consider a few examples of known parametric families for which sampling the pieces will be easy, since the pieces are also from known parametric families. We start with the most important one, the Gamma distribution.

**Definition 3** (Gamma distribution). *A positive, absolutely continuous random variable $X$ is said to be Gamma distributed with shape $\alpha \geq 0$ and scale $s \geq 0$, which we denote $X \sim \mathcal{G}(\alpha, s)$, if it has the density $f(x) = \Gamma(\alpha)^{-1} s^{-\alpha} e^{-\frac{x}{s}} x^{\alpha-1}$, where the normalizing constant, the Gamma function $\Gamma(\alpha)$, is given as $\Gamma(\alpha) = \int_0^\infty y^{\alpha-1} e^{-y} dy$.*

This distribution has mean $\alpha s$ and variance $\alpha s^2$. Its pieces can be computed easily as described in Property 4.

**Property 4** (Pieces of Gamma distributions). *The $\beta$-piece of the $\mathcal{G}(\alpha, s)$ distribution is $\mathcal{G}(\beta \alpha, s)$ distributed.*

*Proof. The distribution has moment generating function given by:*
$$
\begin{aligned}
M(t) = \mathbb{E}(e^{tX}) &= \int_0^\infty e^{tx} f(x) dx \\
&= \Gamma(\alpha)^{-1} s^{-\alpha} \int_0^\infty e^{-x \frac{1-ts}{s}} x^{\alpha-1} dx \\
&= \Gamma(\alpha)^{-1} s^{-\alpha} \int_0^\infty e^{-y} y^{\alpha-1} \left(\frac{1-ts}{s}\right)^{-\alpha} dy \\
&= \Gamma(\alpha)^{-1} s^{-\alpha} \left(\frac{1-ts}{s}\right)^{-\alpha} \int_0^\infty e^{-y} y^{\alpha-1} dy \\
&= (1-ts)^{-\alpha}.
\end{aligned}
$$

*Therefore, Equation (4) gives $M(t)^\beta = (1-ts)^{-\alpha\beta}$ which concludes the proof.* □

From the very structure of its moment generating function $M(t) = (1-ts)^{-\alpha}$, we can already see that the Gamma distribution will be easy to divide. Indeed, the *shape* parameter $\alpha$, which can be any positive real value, is already in the exponent according to Equation (4). This structure also tells us that convolutions of Gamma distributions with same scales are still Gamma distributed (only the *shape* parameter will differ), as Property 4 clearly shows.

There is another distribution class that divides easily, the Gaussian class.

**Definition 5** (Gaussian distribution). *The random variable $X$ is said to be Gaussian with mean $\mu$ and variance $\sigma^2$, which we denote $X \sim \mathcal{N}(\mu, \sigma^2)$, if and only if it has moment generating function $M(t) = e^{\mu t + \frac{\sigma^2 t^2}{2}}$.*

**Property 6** (Pieces of gaussian distribution). *The $\beta$-piece of the $\mathcal{N}(\mu, \sigma^2)$ distribution is $\mathcal{N}(\beta\mu, \beta\sigma^2)$ distributed.*

*Proof. Again, simply compute $M(t)^\beta = e^{\beta\mu t + \frac{\beta\sigma^2 t^2}{2}}$. The same exponent remark applies here.* □

Another interesting example in the discrete case is the Poisson distribution, and most importantly compound Poisson processes.

**Definition 7** (Poisson distribution). *The random variable $N$ has a Poisson distribution with rate $\lambda > 0$, denoted $N \sim \mathcal{P}(\lambda)$, if and only if $\mathbb{P}(N = x) = e^{-\lambda} \frac{\lambda^x}{x!} \mathbb{1}_{x \in \mathbb{N}}$. Its moment generating function is written as $M(t) = e^{\lambda(e^t - 1)}$.*

**Definition 8** (Compound Poisson processes). *Let $Y_1, Y_2, \ldots$ be i.i.d. random variables with distribution $\mathcal{D}$, and let $N \sim \mathcal{P}(\lambda)$. The random variable $X = \sum_{i=1}^{N} Y_i$ is a compound Poisson process, which we denote $X \sim \mathcal{CP}(\lambda, \mathcal{D})$.*

**Property 9** (Pieces of compound Poisson processes). *The $\beta$-piece of the $\mathcal{CP}(\lambda, \mathcal{D})$ distribution is $\mathcal{CP}(\beta\lambda, \mathcal{D})$ distributed.*

*Proof.* We have: $\mathbb{E}(e^{tX}) = \mathbb{E}(\prod_{i=1}^{N} e^{tY_i}) = \mathbb{E}(M_{\mathcal{D}}(t)^N) = e^{\lambda(M_{\mathcal{D}}(t)-1)}$ Hence, $\lambda$ can be seen as a power on the m.g.f. and the same reasoning as before applies. □

Obviously, if $Y_1, Y_2, \ldots$ are all equal to 1, we recover the known fact that the Poisson random variable $N$ is infinitely divisible, with Poisson distributed pieces. Interestingly, the distribution $\mathcal{D}$ of the losses does not matter: since they are i.i.d., thinning the point process is enough to obtain a piece. In fact, the proof applies a little more generally to non-homogeneous Poisson processes: only the independence of increments of the process is needed.

The same proof applies with Negative binomial distributions instead of Poisson distributions: it can be shown that Geometric distributions are infinitely divisible, with negative binomial pieces, which are infinitely divisible with negative binomial pieces themselves.

The next example, the log-Normal distribution, is very important as this distribution is often used in (re)insurance to model losses.

**Definition 10** (Log-Normal distributions). *The random variable $X$ is said to be log-Normally distributed, which we denote $X \sim \mathcal{LN}(\mu, \sigma^2)$, if $\ln X \sim \mathcal{N}(\mu, \sigma^2)$.*

Thorin [10] shows that the log-Normal distribution is infinitely divisible, but does *not* provide the distribution of the pieces, and ends up with the following result:

**Property 11** (Existence of a Thorin measure in the log-Normal case [10]). *The log-Normal distribution is infinitely divisible. Moreover, it can be expressed as a weak limit of Gamma convolutions: there exists a measure $\nu$ such that*

$$M(t) = \exp\left\{-\int_{\mathbb{R}_+} \ln(1 - st)\nu(\partial s)\right\}.$$

Unfortunately, the measure $\nu$, called the Thorin measure, is not explicit in the log-Normal case[1]. Therefore, the log-Normal distributions are not *parametrically* divisible, and obtaining the distributions of the pieces of a log-Normal is a complicated matter.

The same representation can be constructed for Pareto distributions (see [9]), which are also generalized Gamma convolutions with the continuous Thorin measure $\nu$, and therefore not parametrically divisible. Note that the Thorin measure $\nu$ plays the same role as the shape parameter of the Gamma distribution: this representation gives a 'non-parametric' divisibility, where the parameter is a measure $\nu$, not explicitly known, in an infinite-dimensional space, which is not practical.

---

[1] We still know some things about this measure, e.g., that it has infinite total mass. See Bondesson [11] for numerical approximation that allows to plot its density.

In the next section, we propose to approximate this measure $v$ by an atomic measure, or even a simple Dirac (which yields a Gamma distribution). Indeed, a finitely atomic Thorin measure describes a finite convolution of Gamma distribution, which is clearly parametrically divisible with easy-to-sample pieces.

## 3. Approximation schemes

We propose here two of the main parametrically divisible approximations of positive distributions that might be used in additive risk factor model related problems. The first one is a trivial Gamma approximation, while the second one, more involved, is an approximation by a convolution of several Gamma distributions.

### 3-a. The Gamma Approximation

Suppose we have a positive random variable $X$, whose distribution is infinitely divisible but whose pieces are not easy to sample from. Our goal is to approximate the distributions of pieces of $X$ with an approximation that can be easily sampled from. The global strategy is to use an approximating distribution that is parametrically divisible, e.g., a Gamma distribution, to approximate $X$ directly. Indeed, if we have an approximation of $X$ that is Gamma distributed, the pieces can be approximated parametrically through the pieces of this Gamma distribution. Another potential strategy would be to compute the Laplace transform of the distribution and use numerical schemes to sample from an (approximated) Laplace transform, e.g., using a saddle point approximation, as was proposed in [8].

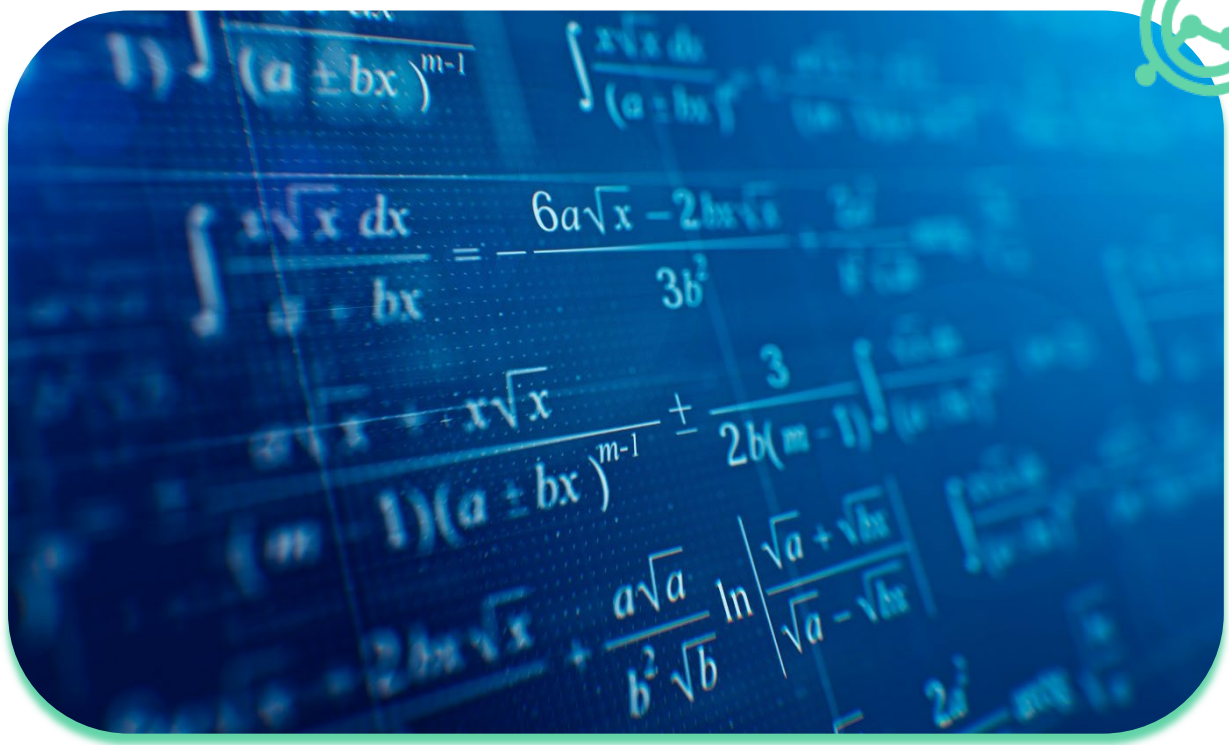

SCOR Paper #43 - Parametric divisibility of stochastic losses    10

We will try to estimate the Gamma distribution by standard maximum likelihood. To do this, recall that the Gamma log-likelihood is given by:

$$l(\alpha, s, x) = -\ln\Gamma(\alpha) - \alpha\ln(s) + (\alpha - 1)\ln(x) - \frac{x}{s}.$$

Hence, for a random variable $X$ such that $\mathbb{E}(X)$ and $\mathbb{E}(\ln X)$ exist, this equation makes it possible to project the distribution of $X$ onto the Gamma class by maximum likelihood.

On the other hand, if either $\mathbb{E}(X)$ or $\mathbb{E}(\ln X)$ does not exist, which might happen with heavy-tailed distributions such as the Pareto distribution (which we will discuss in the next section), and which is quite common in the (re)insurance field, this maximum likelihood approach will not work.

Instead, we propose to match *shifted* moments, i.e., the Taylor coefficients of the moment generating function, but around $-1$ rather than $0$. When $X$ is positive, the moment generating function of $X$, defined by

$$M(t) = \mathbb{E}(e^{-tX}),$$

is analytic in the negative half of the complex plane, and therefore is infinitely derivable in, e.g., $-1$ and has bounded derivatives.

Moreover, if $X \sim \mathcal{G}(\alpha, s)$, the zeroth and first derivative of $M$ are:

$$\begin{aligned}\mu_0 &= M(-1) = (1+s)^{-\alpha} \\ \mu_1 &= M'(-1) = \alpha s (1+s)^{-\alpha-1}\end{aligned}$$

The shifted moments $\mu_0 = \mathbb{E}(e^{-X})$ and $\mu_1 = \mathbb{E}(Xe^{-X})$ could easily be approximated from a sample of the random variable $X$. The set of equations we obtained can then be solved numerically for $\alpha, s$, hence providing an estimator that matches the first two moments of the Escher transform.

Unfortunately, since we only have two parameters $\alpha, s$ we cannot exactly match more than two moments. This is a huge restriction since, as we will illustrate later, the approximation of a log-Normal or a Pareto distribution by a simple Gamma distribution is usually bad. However, by using convolutions of Gamma distributions, it is possible to match more moments and to obtain significantly better approximations, while keeping the parametric divisibility of the estimator. We discuss these possibilities in the next subsection.

### 3-b. Matching more moments with a wider class of distributions

Note that the moment generating function of an independent convolution $Y_1 + \ldots + Y_n$ is the product of the moment generating functions of the convoluted random variables. Therefore, if $Y_1, \ldots, Y_n$ all have $\beta$-pieces given by $Y_1^{(\beta)}, \ldots, Y_n^{(\beta)}$, then $Y_1^{(\beta)} + \ldots + Y_n^{(\beta)}$ is a $\beta$-piece of $Y_1 + \ldots + Y_n$.

Hence, a potentially better estimator than the simple Gamma distribution can be found as a convolution of parametrically divisible distributions. The case that will be of interest to us here is the case of convolutions $n \geq 2$ Gamma distributions. Indeed, we already know from Thorin's work that the log-Normal and Pareto distributions can be seen as limiting cases of convolutions of a finite number of Gamma distributions. Hence, for a certain number of parameters, that is for a certain number of Gammas, we have good reason to hope that the approximation of these distributions by a convolution of Gamma distributions will be almost perfect.

We therefore work with the class of generalized Gamma convolutions.

**Definition 12** (Generalized Gamma convolution [6]). *A random variable $X$ is said to be a generalized Gamma convolution with Thorin measure $v$, denoted $X \sim \mathcal{G}(v)$, if it has a moment generating function*

$$M(t) = \exp\left\{-\int_{\mathbb{R}_+} \ln(1-st)v(\partial s)\right\}.$$

*Moreover, any measure $v$ such that $\int_{[1,\infty)} |\ln s|v(\partial s)l \leq \infty$ and $\int_{(0,1)} s\,v(\partial s) \leq \infty$ generates a valid generalized Gamma convolution.*

We refer to [12–15] for details about the problems involved when estimating such Gamma convolutions and for details about the algorithms that we will use in the next section to estimate them. Note that, of course, finitely atomic measure $v$ fulfills the integration conditions of Definition 12. Moreover, for

$$v = \sum_{i=1}^{n} \alpha_i \, \delta_{s_i},$$

we have

$$M(t) = \prod_{i=1}^{n}(1-s_i t)^{-\alpha_i},$$

which shows that the corresponding $X$ is indeed a convolution of $n$ Gamma distributions, with respective parameters $(\alpha_1, s_1), \ldots, (\alpha_n, s_n)$.

Finally, the shape of $M(t)$ show that the parametric division of a convolution of Gamma distributions will be quite easy, using the same exponent trick as before.

**Property 13** (Pieces of generalized Gamma convolutions). *The $\beta$-piece of a generalized Gamma convolution $X \sim \mathcal{G}(v)$ is also a generalized Gamma convolution with Thorin measure $\beta v$.*

## 4. Illustrative examples

We propose a quick illustration of the potential applications of an additive risk factor model as defined by Equations (1), (2) and (3), to be able to show the dependency that can arise from such a model, and to assess the performance of the proposed Gamma approximation. We start with the description of a simple division of a distribution, and we end with the sampling from a genuine additive risk factor model using previously discussed parametric approximations.

## 4-a. Division of the Pareto distribution

The Pareto distribution is used in (re)insurance loss modeling due to its heavy tail. It also has the advantage of being an extreme value distribution, and is therefore well suited for losses that are defined as the maximum of other random variables.

**Definition 14** (Pareto distribution). *We say that the random variable $X$ with support $\mathbb{R}_+$ is Pareto distributed with shape $\alpha \in \mathbb{R}_+$ if it has the density*

$$f(x) = \alpha(x+1)^{-\alpha-1}.$$

Recall that the shape parameter influences the integrability of the distribution: the variance is finite if and only if $\alpha > 2$, the expectation is finite if and only if $\alpha > 1$, and the density is square integrable if and only if $\alpha > \frac{1}{2}$. Indeed, the smaller the value of $\alpha$, the heavier the tail of the distribution.

The Pareto distribution is infinitely divisible, as shown by Thorin in [9]. Unfortunately, there is no parametric form for the $\beta$-piece of a Pareto distribution. On the other hand, Thorin shows that the Pareto distributions belong to the wider class of generalized Gamma convolutions, and can therefore be approximated arbitrarily well by a finite convolution of Gamma distributions, which we can easily take pieces of to sample from.

Our goal is to divide the Pareto distribution into pieces. For that, we propose two approximations. The easiest approximation that we will discuss is the Gamma approximation. We will also consider the convolution of $n$ Gammas, $n \geq 2$. According to the previous section, both these constructions are parametrically divisible.

As our example, we consider the random variable $X$ as having a Pareto distribution with parameter $\alpha = \frac{3}{4}$. This distribution has no expectation or variance, and is not parametrically divisible. However, it can be approximated by a parametrically divisible model, either a simple Gamma distribution, or a convolution of, for example, $n = 20$ Gamma distributions. The first approximation is obtained using a maximum likelihood principle on a large sample, and the second one is obtained according to the method from [12].

The two following graphs describe the quality of both these parametrically divisible approximations of $X$.

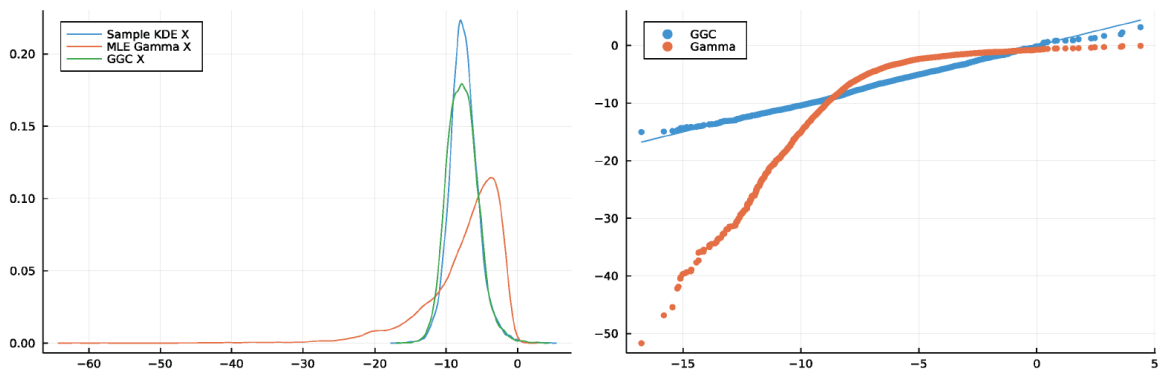

*Figure 1: (a) Density plots of the Pareto density (via a sample and a kernel estimation), its projection into the Gamma class and into the GGC class with $n = 20$. The x-axis is*

*log-scaled. (b) Quantile-Quantile plots of the Gamma approximation and the GGC approximation (convolution of $n = 20$ Gammas) against the true Pareto distribution.*

We see from the density plot on Figure 1(a) and the quantile-quantile plot on Figure 1(b) that the quality of the convolution of $n = 20$ Gamma is almost perfect, reproducing the Pareto distribution completely up to very high quantiles. Conversely, the Gamma approximation struggles. Note that a good approximation of the distribution itself will induce well-approximated pieces.

To better understand the difference between the two approximations and the impact it will have on a dependence structure built through a risk-factor model with Pareto margins, we propose in the next section to sample from a multivariate model with a dependence structure given as an additive risk factor model, using this pareto distribution as one of the marginals.

### 4-b. Approximate sampling from a risk factor model

Consider that we furthermore observe a second random variable $Y$, distributed as a log-Normal distribution with log-mean $\mu = 0$ and log-standard deviation $\sigma = 2$. This distribution also belongs to the Thorin class of generalized Gamma convolutions (see [10]), and therefore we can estimate it as a Gamma distribution or as a convolution of $n = 20$ Gamma distributions in the same way we did for $X$. The result of these approximations is given in Figure 2.

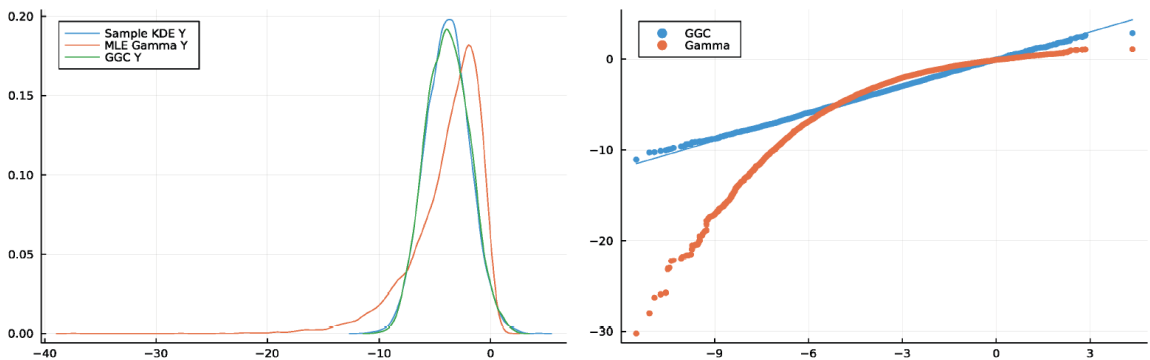

*Figure 2: (a) Density plots of the log-Normal density (via a sample and a kernel estimation), its projection into the Gamma class and into the GGC class with $n = 20$. The x-axis is log-scaled. (b) Quantile-Quantile plots of the Gamma approximation and the GGC approximation (convolution of $n = 20$ Gammas) against the true log-Normal distribution.*

We observe, as for the Pareto distribution, that the Gamma approximation is not satisfactory, while the convolution of 20 Gammas is essentially perfect at this level of detail (the quantile-quantile plot has 10 000 points, hence our approximation seems accurate up to the $99.99\%$ quantile). Since the Gamma convolution approximations of the log-Normal and Pareto distributions are almost perfect, we will benchmark the results obtained by the Gamma approximation against them.

We consider a dependency between $X$ and $Y$ that is given by a risk factor model as follows: we consider $X_{0.2}, X_{0.8}, Y_{0.2}$ and $Y_{0.8}$ the 20% and 80%-pieces of $X$ and $Y$, and we suppose that $X_{0.8}$ and $Y_{0.2}$ are comonotone. Our goal is to approximately sample from the obtained random vector $(X, Y)$. More precisely, the model has the following structure:

$$\begin{aligned} X &= X_{0.2} + X_{0.8} \\ Y &= Y_{0.8} + Y_{0.2} \end{aligned}$$
$X_{0.2}, Y_{0.8}$ and $(X_{0.8}, Y_{0.2})$ are mutually independent
$X_{0.8}$ and $Y_{0.2}$ are comonotone

Note that this specifies completely the distribution of the random vector $(X, Y)$.

Since the Gamma approximation and the GGC approximation are parametrically divisible, constructing approximated distributions for their pieces is quite easy. Then, after sampling the four pieces, we reorder the samples so that $X_{0.8}$ and $Y_{0.2}$ samples have the same ranks, and we sum back to obtain $X$ and $Y$. The obtained dependence structures, on a copula-scale, can be seen in Figure 3.

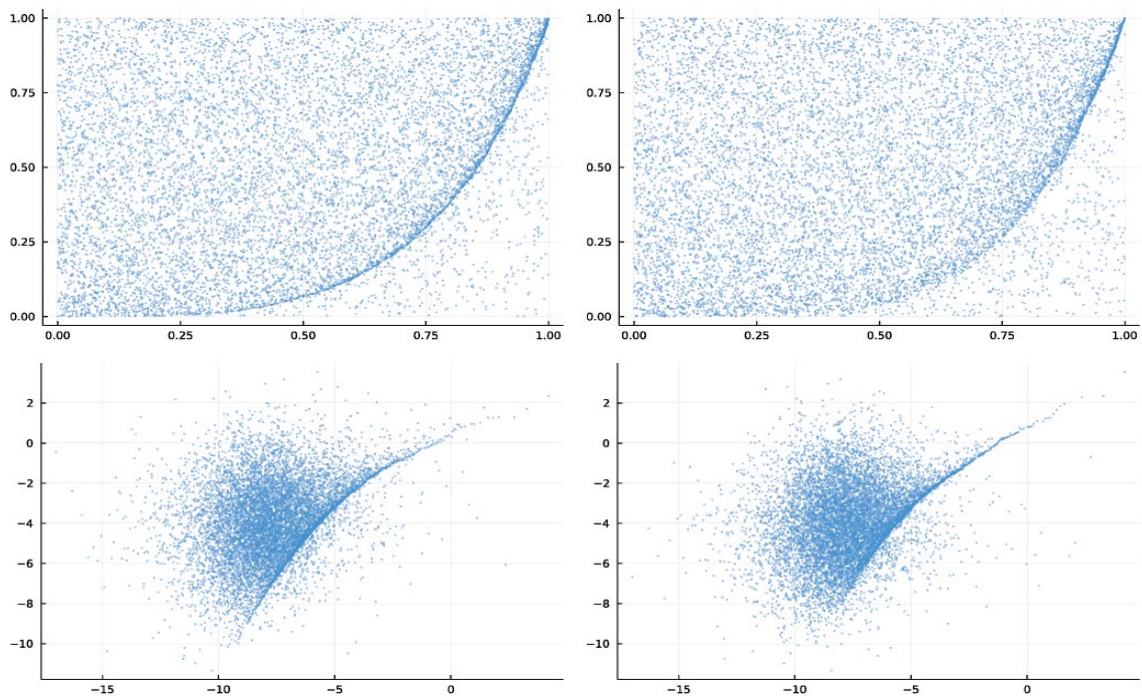

Figure 3: On the left, a sample from the Gamma approximation, while on the right is a sample from the convolution of Gamma approximation given as a benchmark. The top row presents the samples on the copula scale (pseudo-observations), while the bottom row presents them on the log scale.

We note that the dependence structure produced by the Gamma approximation is clearly not the right one. From these obtained empirical ranks, we use the quantile functions of the initial Pareto and log-Normal distributions to produce a sample, whose scatter plot is visible on the bottom row of Figure 3. Surprisingly, we observe very little difference between the scatter plots on the marginals scale. Indeed, the dependence structure error is quite small. This is moreover confirmed by a kernel density estimation on $X + Y$, shown in Figure 4.

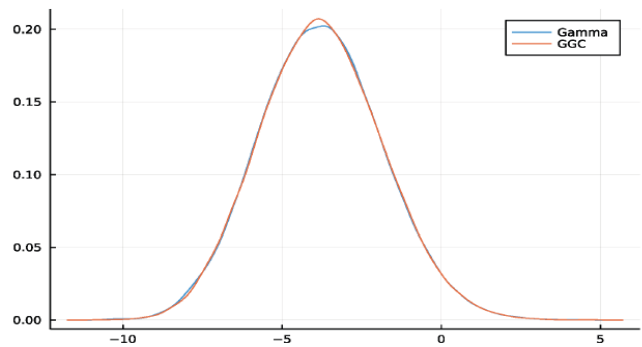

*Figure 4: Density of the sum $X + Y$ obtained from both approximations.*

Figure 4 shows that, at least in this particular example, the use of a bad approximation to produce the pieces does not affect the produced dependence structure enough to alter the distribution of the total sum: the error seems negligible with respect to the sampling noise.

If we consider the convolution of Gamma approximation to be the ground truth, which is credible because it approximates our distributions very well (see [12] for detailed experiments), we can conclude that, at least in this simple case, the dependence structure obtained from the Gamma approximation is close enough to the truth to allow the production of meaningful samples of the random vector when paired with the true marginals.

## 5. Conclusion

A risk factor model is a useful dependency structure in (re)insurance internal modeling because its interpretability facilitates the calibration of the dependence structures in those cases lacking data. However, sampling such a model, which is often necessary to derive explicit quantities useful to capital management, such as the upper quantiles of the total risk, can be really difficult.

This difficulty comes almost directly from the challenges involved in sampling a given piece of a given divisible distribution. Therefore, some kind of approximation must be used to allow for practical implementations. Such approximation will perform correctly if it is parametrically divisible, that is if its pieces have a known parametric distribution, or at least can easily be sampled from.

From the two approximations we proposed here, the most precise is of course the expansion into the Thorin class of generalized Gamma convolutions. The Thorin class does not contain all infinitely divisible distributions, but it does contain a lot of common shapes, including the heavy tailed distributions that are so useful in insurance applications. However, this expansion can be numerically challenging (as observed back in 2019 by [16], and developed largely in [12-14]). On the other hand, the simpler Gamma approximation is easy to implement and still provides a dependence structure that is quite close to the one fixed in the original risk factor model. The closer our distributions are to Gamma distributions, the closer the approximated dependence structure will be to the truth.

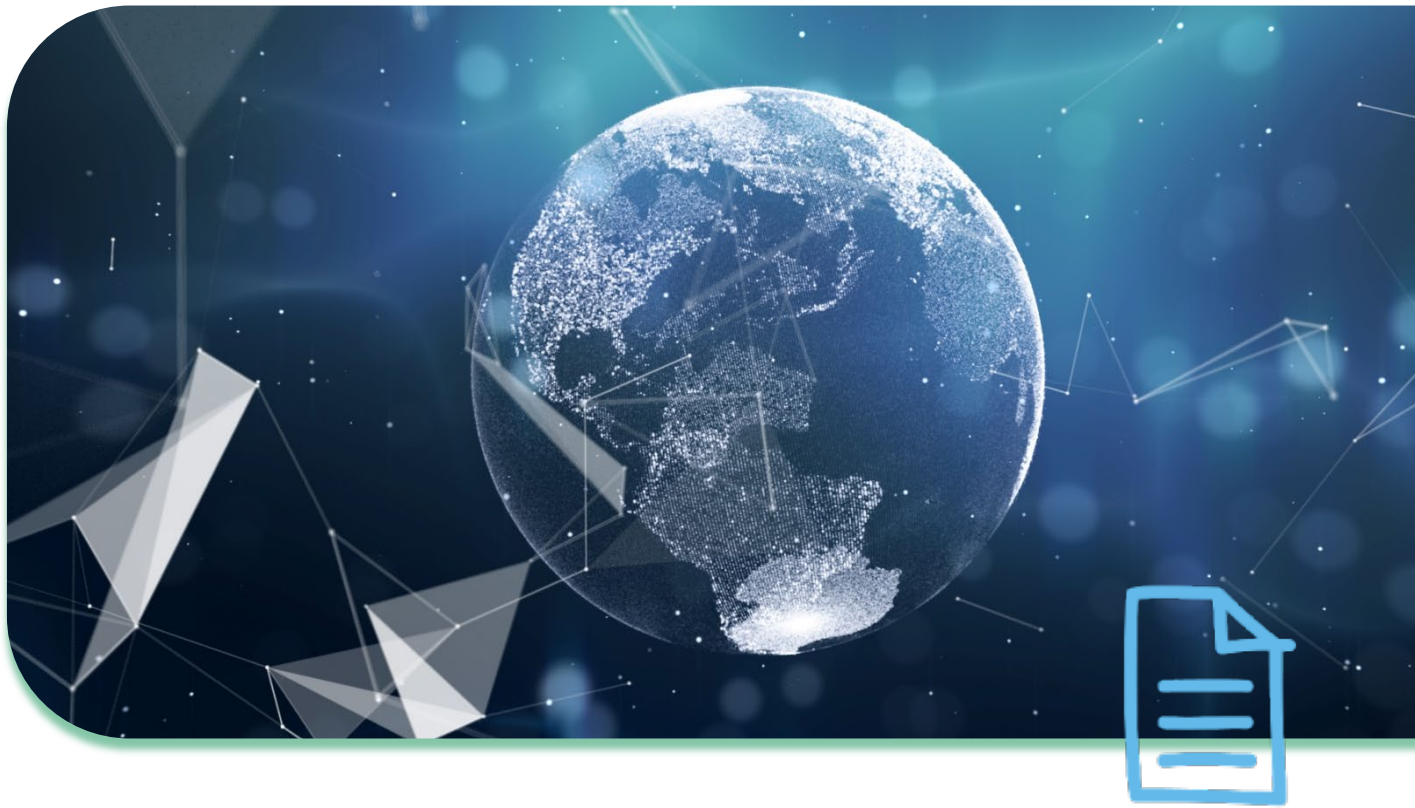